\documentclass{article}

    \usepackage[preprint]{neurips_2025}

\usepackage[utf8]{inputenc} 
\usepackage[T1]{fontenc}    
\usepackage{hyperref}       
\usepackage{url}            
\usepackage{amsfonts}       
\usepackage{nicefrac}       
\usepackage{microtype}      
\usepackage{makecell}
\usepackage{algorithm}
\usepackage{algpseudocode}
\usepackage{booktabs}
\usepackage{multirow} 
\usepackage{array} 

\usepackage{amsmath,amssymb,amsfonts,amsthm,bbm}
\usepackage{mathtools}   
\usepackage{stmaryrd} 

\usepackage{enumitem}   
\usepackage{romannum}

\usepackage[titletoc,title]{appendix}  

\usepackage[table,xcdraw,dvipsnames]{xcolor}
\usepackage{hyperref}
\newcommand\myshade{85}
\colorlet{mylinkcolor}{YellowOrange}
\colorlet{mycitecolor}{Aquamarine}
\colorlet{myurlcolor}{violet}

\hypersetup{
  linkcolor  = mylinkcolor!\myshade!black,
  citecolor  = mycitecolor!\myshade!black,
  urlcolor   = myurlcolor!\myshade!black,
  colorlinks = true,
}

\usepackage{caption}
\usepackage{subcaption}

\usepackage[normalem]{ulem}
\usepackage{setspace}
\usepackage{graphicx}
\usepackage{comment}
\graphicspath{{figures/}}
\title{Joint Optimization of Service Routing and Scheduling in Home Health Care}

\author{%
Yi Zhang \\
Columbia University \\
New York, NY 10027 \\
\texttt{yz5195@columbia.edu} \\
\AND
Zhenzhen Zhang \\
Tongji University \\
Shanghai, 20092 \\
\texttt{zhenzhenzhang@tongji.edu.cn} \\
}

\begin{document}
\pagenumbering{arabic}

\maketitle

\begin{abstract}
The growing aging population has significantly increased demand for efficient home health care (HHC) services. This study introduces a \textit{Vehicle Routing and Appointment Scheduling Problem} (VRASP) to simultaneously optimize caregiver routes and appointment times, minimizing costs while improving service quality. We first develop a deterministic VRASP model and then extend it to a stochastic version using sample average approximation to account for travel and service time uncertainty. A tailored Variable Neighborhood Search (VNS) heuristic is proposed, combining regret-based insertion and Tabu Search to efficiently solve both problem variants. Computational experiments show that the stochastic model outperforms the deterministic approach, while VNS achieves near-optimal solutions for small instances and demonstrates superior scalability for larger problems compared to CPLEX. This work provides HHC providers with a practical decision-making tool to enhance operational efficiency under uncertainty.

\end{abstract}

\section{Introduction}\label{sec:intro}
According to research by the United Nations Department of Economic and Social Affairs, population aging has become a defining global trend. In China, by the end of 2023, the population aged 60 and above reached 297 million (21.1\% of total population), with 209 million (14.9\%)~\cite{1} aged 65 and above. Globally, the population aged 65+ was 761 million in 2021 and is projected to reach 1.6 billion by 2050~\cite{2}. This demographic shift has significantly increased demand for Home Health Care (HHC) services, which provide medical tests, wound care, and nursing visits at patients' homes~\cite{3,4}.

Traditional HHC systems face operational challenges where either healthcare workers or customers experience unnecessary waiting due to rigid appointment scheduling~\cite{5}. This inefficiency stems from allowing customers to unilaterally determine service times, which often leads to suboptimal routing plans. A more effective approach would enable medical institutions to jointly optimize service schedules and routing paths through negotiation with clients - a model already adopted by some HHC providers~\cite{6,7,8,9}.

The problem is further complicated by real-world uncertainties. Healthcare workers typically serve multiple geographically dispersed clients per trip, where travel times fluctuate due to traffic conditions and service durations vary unexpectedly. These uncertainties frequently cause either worker idle time or client waiting periods. 

To address these challenges, this study proposes a Vehicle Routing and Appointment Scheduling Problem (VRASP) that simultaneously optimizes service scheduling and routing while accounting for travel and service time uncertainties. Our joint optimization approach aims to improve resource utilization and service punctuality in HHC operations.

\section{Related Work}\label{sec:related}

\subsection{Joint Optimization in Vehicle Routing Problems}

The Vehicle Routing Problem (VRP), first introduced by Dantzig and Ramser in 1959, seeks optimal delivery routes from a depot to multiple customers under capacity constraints~\cite{10,11,12}. Its variants include the Vehicle Routing Problem with Time Windows (VRPTW), which incorporates either hard time windows (mandatory arrival times) or soft time windows (penalized deviations)~\cite{13}.

Joint optimization simultaneously optimizes multiple interdependent decision variables under various constraints, offering comprehensive system insights and revealing variable dependencies that lead to more scientific and effective decisions~\cite{24}. Recent research has identified three main types of joint optimization in vehicle routing problems~\cite{23}: (i) Inventory-Routing Problems (IRP): Simultaneously optimize delivery quantities and transportation routes to prevent stockouts; (ii) Location-Routing Problems (LRP): Jointly determine facility locations and distribution routes, with applications ranging from waste collection to electric vehicle charging station planning; (iii) Routing-Driver Scheduling: Integrate driver assignment and shift scheduling with route optimization, though this remains less studied compared to other variants.

Bell et al.~\cite{25} first proposed this class of problems when studying delivery scheduling and transportation route decisions for chemical companies. Subsequently, Federgruen and Zipkin~\cite{26} investigated how to determine both resource allocation at multiple locations and fleet delivery routes under uncertain customer demand. The Location Routing Problem (LRP) aims to simultaneously determine warehouse locations and routes between warehouses that serve geographically dispersed customers. For instance, Fischer et al.~\cite{27} developed a two-phase heuristic algorithm for waste collection scenarios, where the first phase selects optimal collection points based on customer location distribution, and the second phase determines minimum-cost delivery routes from collection points to customers. Yao et al.~\cite{28} studied electric vehicle scheduling problems, considering charging requirements during trips, which requires simultaneous optimization of routing and charging station selection - a variant of LRP. The Routing and Driver Scheduling Problem, which integrates driver assignment with route planning, represents an important research trend in transportation management, though literature remains relatively sparse compared to other variants~\cite{23}. Geol~\cite{29} examined work-rest schedules and routing for truck drivers, proposing a large neighborhood search algorithm. Wen et al.~\cite{30} addressed vehicle routing and driver scheduling in supermarket fresh meat supply chains using a multilevel variable neighborhood search heuristic.

While joint optimization has been extensively studied in VRPTW variants, no existing work has investigated the simultaneous optimization of service start times and vehicle routing - the novel contribution of our research. Our approach addresses this gap by treating service scheduling and routing as interdependent decision variables, enabling more efficient resource allocation in home healthcare logistics. The proposed methodology demonstrates particular advantages in scenarios requiring precise time coordination between multiple service providers and customers.

\subsection{Uncertainty in Vehicle Routing Problems}

Traditional deterministic vehicle routing models assume all parameters (travel times, service durations, etc.) are known precisely, but real-world applications inevitably face various uncertainties. Research on uncertain VRPTW can be categorized into four types: (i) uncertain customer demand, (ii) uncertain customer presence, (iii) uncertain service times and travel times. While the first two categories have been extensively studied~\cite{33,34,35}, less attention has been paid to service and travel time uncertainties - the focus of this work. Two primary approaches exist for modeling uncertainty: robust optimization and stochastic programming.

\paragraph{Robust Optimization} In robust VRPTW, uncertain parameters are defined by uncertainty sets. Cho et al.~\cite{36} investigated the TSP with time windows problem and addressed the potential over-conservatism in robust optimization by introducing a conservatism degree parameter to their model. This approach maintains computational efficiency while providing control over solution robustness. Montemanni et al.~\cite{37} modeled a similar problem as a robust shortest path problem on interval graphs, proposing a branch-and-bound algorithm for its solution. Agra et al.~\cite{38} studied VRPTW with uncertain travel times using uncertainty set modeling, developing feasible solutions for all possible travel time scenarios through combined branch-and-cut and column-and-row generation algorithms. For the same problem, Wu et al.~\cite{39} introduced a novel robustness metric and validated its superiority over classical measures (e.g., worst-case, best-case, and minimax-regret criteria) through extensive computational experiments. Nasri et al.~\cite{40} simultaneously considered travel time and service time uncertainties, employing parallel processing to generate scenario sets and an adaptive large neighborhood search heuristic to identify optimal solutions for each scenario.

Jaillet et al.~\cite{41} proposed the Requirements Violation Index (RVI) based on Aumann and Serrano's~\cite{42} Riskiness Index (RI) as a more comprehensive metric for evaluating tardiness impacts, incorporating both probability and magnitude of delays. However, this method is only applicable when arrival times can be expressed as affine functions of independent factors - a condition violated in classical VRPTW where arrival times are piecewise linear functions of travel times. To address this limitation, Zhang et al.~\cite{43} developed the Essential Riskiness Index (ERI) for VRPTW with empirically distributed travel times, though this approach only guarantees optimality for small-scale instances. Further advancing this research stream, Zhang et al.~\cite{44} introduced the Service Fulfillment Risk Index (SRI) and formulated a Distributionally Robust Optimization (DRO) model using Wasserstein distance-based ambiguity sets to prevent overfitting to empirical travel time distributions. Their numerical experiments demonstrated significant improvements in both computational efficiency and solution quality.

\paragraph{Stochastic Optimization} For stochastic VRPTW, where travel and service times follow known continuous or discrete distributions~(e.g., uniform, normal), Laporte et al.~\cite{45} established a chance-constrained model minimizing operational costs while restricting overtime probability. They proved that normality assumptions enable reformulation of chance constraints as deterministic equivalents. Kenyon and Morton~\cite{46} proposed two branch-and-cut algorithms: exact solutions for small discrete sample spaces and sampling-based heuristics for larger instances, with objectives focusing respectively on minimizing expected return time and maximizing on-time return probability. Notably, these studies primarily optimized travel costs without considering penalties for time window violations.

Taniguchi et al.~\cite{47} developed a Vehicle Routing and Scheduling Problem with Time Windows~(VRPTW-P) incorporating both operational costs and early/tardy penalties. Using dynamic traffic simulation to estimate travel time distributions, their model imposes proportional penalties for early arrivals (provider waiting costs) and late arrivals~(client delay costs). Ando and Taniguchi~\cite{48} extended this work with enhanced travel time uncertainty analysis, estimating operational costs through mean travel times and computing penalties via arrival time distributions. Their genetic algorithm implementation revealed significant impacts of travel time uncertainty on delay costs. Russell and Urban~\cite{49} studied stochastic VRPTW with soft time windows using offset gamma distributions for travel times and Taguchi loss functions for penalty calculation, solved via tabu search. Branda~\cite{50} investigated hard time window variants through mixed-integer stochastic programming with sample approximation techniques. Expósito et al.~\cite{51} integrated random travel and service times in a QoS-oriented model solved by a GRASP-VNS hybrid metaheuristic. Ehmke et al.~\cite{52} proposed a chance-constrained model limiting individual time window violation probabilities through service start time distribution analysis, enabling efficient large-scale solutions when combined with existing VRPTW algorithms.

\subsection{Sample Average Approximation}
The Sample Average Approximation (SAA) method is a technique that uses Monte Carlo simulation to solve stochastic optimization problems. In this method, the expected objective function of a stochastic problem is approximated by the sample average estimate from random samples, thereby transforming the stochastic programming problem into a deterministic one. The process can be further repeated with different samples to obtain candidate solutions and statistical estimates of their optimality gaps~\cite{53}. In the interior sampling method, the samples drawn continuously change as the optimization process progresses. In contrast, the SAA method is an exterior sampling method, where N samples are generated according to probability distribution P, and then the deterministic optimization problem specified by the drawn samples is solved. This process can be repeated through multiple samplings and solutions of the specified deterministic optimization problem.

The application of the SAA method in stochastic programming problems began in the 1990s, with significant research conducted by Geyer and Thompson~\cite{54} and Rubinstein and Shapiro~\cite{55}. Subsequently, Plambeck~\cite{56} and Kleywegt et al.~\cite{57} further studied the SAA method in stochastic linear programming and discrete stochastic programming, respectively. Beyond theoretical research, some scholars have applied the SAA method to practical problems. For example, Ng et al.~\cite{58} used the SAA method to solve the stochastic resource optimization problem in wireless-powered hybrid coded edge computing networks. Rath and Rajaram~\cite{59} studied hospital staff scheduling, modeling the problem as a two-stage stochastic integer dynamic programming model and solving it using the SAA method. In the context of vehicle routing problems (VRP) and related issues, the SAA method has also been employed by some researchers. For instance, Verweij et al.~\cite{60} effectively solved the stochastic TSPD problem by combining Benders decomposition and branch-and-cut algorithms. Tas et al.~\cite{61} investigated the VRPTW with soft time windows where travel times follow a gamma distribution, designing a branch-and-bound algorithm to minimize operational costs and penalties for time window violations.

\section{Deterministic VRASP Model}\label{sec:prob-setting}
\subsection{Problem Description}
Consider a healthcare center with multiple caregivers who need to visit clients during working hours. Each caregiver departs from the center, follows an assigned route to serve clients sequentially, and returns to the center. Overtime costs are incurred if a caregiver returns after the designated working hours.

Clients submit service requests in advance, and the healthcare center must determine both the visiting routes and scheduled service start times for each client. If a caregiver arrives earlier than the scheduled time, they must wait; late arrivals immediately begin service but incur penalty costs for client waiting time. The model aims to optimize caregiver routing and service scheduling to minimize the sum of fixed operational costs, travel costs, late arrival penalties and overtime costs.

\begin{table}[htbp]
\renewcommand{\arraystretch}{1.2}
\setlength{\tabcolsep}{12pt}
\centering
\caption{Model Parameters and Variables}
\begin{tabular}{ll}
\hline
Symbol & Description \\ \hline
\multicolumn{2}{l}{\textbf{Parameters}} \\ \hline
$N$ & Set of all clients, $N = \{1, 2, 3, ..., |N|\}$ \\ 
$0, |N| + 1$ & Origin and destination (both represent the healthcare center) \\ 
$V$ & Set of all nodes, $V = N \cup \{0\} \cup \{|N| + 1\}$ \\ 
$V_1$ & Set of all clients and origin, $V_1 = N \cup \{0\}$ \\ 
$V_2$ & Set of all clients and destination, $V_2 = N \cup \{|N| + 1\}$ \\ 
$K$ & Set of all caregivers, $K = \{1, 2, 3, ..., |K|\}$ \\ 
$t_{ij}$ & Travel time from client $i$ to client $j$ \\ 
$ts_i$ & Service time at client $i$ \\ 
$c^f$ & Fixed scheduling cost per caregiver \\ 
$c^o$ & Unit overtime cost for exceeding working hours \\ 
$c_{ij}^r$ & Travel cost from client $i$ to client $j$ \\ 
$c^t$ & Unit tardiness cost for late arrival \\ 
$[0, L]$ & Working time range for caregivers \\ \hline
\multicolumn{2}{l}{\textbf{Decision Variables}} \\ \hline
$s_i$ & Scheduled service start time for client $i$ \\ 
$x_{ijk}$ & Binary variable indicating if caregiver $k$ travels from client $i$ to $j$ \\ \hline
\multicolumn{2}{l}{\textbf{Auxiliary Variables}} \\ \hline
$a_i$ & Actual service start time at client $i$ \\ 
$w_i$ & Waiting time at client $i$ \\ 
$o_k$ & Overtime duration for caregiver $k$ \\ \hline
\end{tabular}
\end{table}

\subsection{Deterministic VRASP Formulation}
In deterministic VRASP, where travel times $t_{ij}$ and service times $ts_i$ are known constants, we formulate model $P_0$ as follows:
\begin{align}
\min \quad & \sum_{j \in N} \sum_{k \in K} c^f x_{0jk} + \sum_{k \in K} \sum_{i \in V_1} \sum_{j \in V_2} c_{ij}x_{ijk} + \sum_{i \in V} c_i^t w_i + \sum_{k \in K} c^o o_k \label{eqn:1} \\
\text{s.t.} \quad & \sum_{k \in K} x_{ijk} = 1, \quad \forall i \in N, \label{eqn:2} \\
& \sum_{j \in V_2} x_{0jk} = 1, \quad \forall k \in K, \label{eqn:3} \\
& \sum_{i \in V_1} x_{i(|N|+1)k} = 1, \quad \forall k \in K, \label{eqn:4} \\
& \sum_{j \in V_2} x_{ijk} = \sum_{j \in V_1} x_{jik}, \quad \forall i \in N, \forall k \in K, \label{eqn:5} \\
& w_j = \max(s_j - a_j, 0), \quad \forall j \in V, \label{eqn:6} \\
& \left( \max(a_i, s_i) + ts_i + t_{ij} \right) - a_j \leq M(1 - x_{ijk}), \quad \forall i \in V, \forall j \in V, \forall k \in K, \label{eqn:7} \\
& s_i + ts_i + t_{ij} - s_j \leq M(1 - x_{ijk}), \quad \forall i \in V, \forall j \in V, \forall k \in K, \label{eqn:8} \\
& s_i + ts_i + t_{i(|N|+1)} - L - o_k \leq M(1 - x_{i(|N|+1)k}), \quad \forall i \in V, \forall k \in K, \label{eqn:9} \\
& x_{ijk} \in \{0, 1\}, s_i \geq 0, \quad \forall i \in V, \forall j \in V, \forall k \in K, \label{eqn:10} \\
& a_i, w_i, o_k \geq 0, \quad \forall i \in V, \forall k \in K, \label{eqn:11}
\end{align}

Constraint~\eqref{eqn:1} minimizes the total cost including fixed, travel, tardiness, and overtime costs. Constraint~\eqref{eqn:2} ensures each client is served by exactly one caregiver. Constraints~\eqref{eqn:3}-\eqref{eqn:5} define route continuity from origin through clients to destination. Constraint~\eqref{eqn:6} calculates client waiting time. Constraints~\eqref{eqn:7}-\eqref{eqn:8} establish temporal relationships between consecutive visits. Constraint~\eqref{eqn:9} computes overtime duration. Constraints~\eqref{eqn:10}-\eqref{eqn:11} define variable domains.

\section{Stochastic VRASP Model}
\subsection{Stochastic VRASP Formulation}
Real-world scenarios often involve significant uncertainties. Caregivers typically serve multiple geographically dispersed clients in each trip, where travel times may vary due to traffic congestion, and service durations may fluctuate due to unpredictable factors during service delivery. These uncertainties can lead to either caregiver idle time (if arrivals are early) or client waiting time (if arrivals are delayed). Incorporating travel and service time uncertainties into the decision-making process enables better risk management and more practical solutions, improving both caregiver utilization and service punctuality. Building upon the deterministic model, we develop a stochastic VRASP model using the Sample Average Approximation (SAA) method.

\paragraph{Sample Average Approximation} In stochastic VRASP, both travel times $t_{ij}$ and service times $ts_i$ are random parameters. We employ the SAA method with Monte Carlo sampling to transform the stochastic problem into a deterministic equivalent with finite scenarios, thereby reducing computational complexity. The notation for sampled parameters is summarized in Table~\ref{table:saa}.

\begin{table}[h]\label{table:saa}
\renewcommand{\arraystretch}{1.2} 
\setlength{\tabcolsep}{12pt} 
\centering
\caption{SAA Notation}
\begin{tabular}{cl}
\hline
Symbol & Description \\ \hline
$M$ & Number of samples \\ 
$t_{ij}^m$ & Realized travel time from client $i$ to $j$ in sample $m$ \\ 
$ts_i^m$ & Realized service time at client $i$ in sample $m$ \\ 
$a_i^m$ & Realized service start time at client $i$ in sample $m$ \\ 
$w_i^m$ & Realized waiting time at client $i$ in sample $m$ \\ 
$o_k^m$ & Realized overtime duration for caregiver $k$ in sample $m$ \\ \hline
\end{tabular}
\end{table}

Through SAA, the original problem\eqref{eqn:1}-\eqref{eqn:11} transforms into stochastic problem~\eqref{eqn:12}-\eqref{eqn:17}, denoted as $P_1$:

\begin{align}
\min \quad & \sum_{j \in N} \sum_{k \in K} c^f x_{0jk} + \sum_{k \in K} \sum_{i \in V_1} \sum_{j \in V_2} c_{ij}x_{ijk} + \frac{1}{|M|} \sum_{m=1}^{M} \left( \sum_{i \in V} c_i^t w_i^m + \sum_{k \in K} c^o o_k^m \right) \label{eqn:12} \\
\text{s.t.} \quad & \text{\eqref{eqn:2}-\eqref{eqn:5}, \eqref{eqn:10}} \nonumber \\
& w_j^m = \max(s_j - a_j^m, 0), \quad \forall j \in V, \label{eqn:13} \\
& \left( \max(a_i^m, s_i) + ts_i^m + t_{ij}^m \right) - a_j^m \leq M(1 - x_{ijk}), \quad \forall i,j \in V, \forall k \in K, \label{eqn:14} \\
& s_i + ts_i^m + t_{ij}^m - s_j \leq M(1 - x_{ijk}), \quad \forall i,j \in V, \forall k \in K, \label{eqn:15} \\
& s_i + ts_i^m + t_{i(|N|+1)}^m - L - o_k^m \leq M(1 - x_{i(|N|+1)k}), \quad \forall i \in V, \forall k \in K \label{eqn:16} \\
& a_j^m, w_j^m, o_k^m \geq 0, \quad \forall i \in V, \forall k \in K, \label{eqn:17}
\end{align}

\subsection{SAA Computational Procedure}
Generate $|Q|$ independent sample sets $M_1, M_2, ..., M_{|Q|}$, each containing $|M|$ i.i.d. samples:

\[ M_q = \{ \omega_q^1, \omega_q^2, ..., \omega_q^{|M|} \}, \quad q = 1, 2, ..., |Q| \]

For each sample set $M_q$, solve model $P_1$ to obtain (i) optimal schedules $s_M^1, s_M^2, ..., s_M^Q$, (ii) staffing solutions $x_M^1, x_M^2, ..., x_M^Q$, and (iii) minimum costs $c_M^1, c_M^2, ..., c_M^Q$.

The mean optimal value provides a lower bound (LB) estimate:
\begin{equation*}
    \bar{c}_M = \frac{1}{|Q|} \sum_{q=1}^{Q} c_M^q , \quad  E(\bar{c}_M) \leq c^{*}
\end{equation*}

For any feasible solution $\hat{x} \in X$, the upper bound (UB) is estimated by:

\begin{align*}
\hat{c}_{M'}(\hat{x}) &= \sum_{j \in N} \sum_{k \in K} c^f x_{0jk} + \sum_{k \in K} \sum_{i \in V_1} \sum_{j \in V_2} c_{ij}x_{ijk} + \frac{1}{|M'|} \sum_{m=1}^{|M'|} \left( \sum_{i \in V} c_i^t w_i^m + \sum_{k \in K} c^o o_k^m \right) 
\end{align*}

where $|M'| \gg |M|$ ensures $\hat{c}_{M'}(\hat{x})$ is an unbiased UB estimator with $E(\hat{c}_{M'}(\hat{x})) \geq c^{*}$.

The variances are calculated as:

\begin{align*}
\hat{\sigma}_{\bar{c}_M}^2 &= \frac{1}{|Q|(|Q|-1)} \sum_{q=1}^{Q} (c_M^q - \bar{c}_M)^2 \\
\hat{\sigma}_{\hat{c}_{M'}(\hat{x})}^2 &= \frac{1}{|M'|(|M'|-1)} \sum_{m=1}^{M} \left( \sum_{j,k} c_j x_{0jk} + \sum_{i,j,k} c_{ij} x_{ijk} + \sum_i c_i w_i^m + \sum_k c^o o_k^m - \hat{c}_{M'}(\hat{x}) \right)^2
\end{align*}

From the $|Q|$ candidate solutions $x_M^1, ..., x_M^Q$, select the optimal estimate:

\[ \hat{x}^* \in \arg \min \{ \hat{c}_{M'}(\hat{x}) : \hat{x} \in \{ x_M^1, ..., x_M^Q \} \} \]

Evaluate solution quality using the optimality gap:

\begin{align*}
\text{gap} &= \hat{c}_{M'}(\hat{x}^*) - \bar{c}_M \\
\hat{\sigma}_{\text{gap}}^2 &= \hat{\sigma}_{\bar{c}_M}^2 + \hat{\sigma}_{\hat{c}_{M'}(\hat{x})}^2
\end{align*}

This evaluation framework was originally proposed by Norkin et al.~\cite{63} and further developed by Mak et al.~\cite{64}, with convergence properties analyzed by Kleywegt et al~\cite{65}.

\section{Variable Neighborhood Search Algorithm for VRASP}

\subsection{An Overview of Our VRASP ALlgorithm}
This study designs a Variable Neighborhood Search (VNS) algorithm to solve the VRASP. Algorithm 1 presents the overall framework of the VNS algorithm, where $N = \{N_1, N_2, N_3\}$ represents the neighborhood set, and $N'$ denotes the available neighborhoods during the search process. 

The VNS starts with an initial solution generated by the Clarke-Wright (C-W) algorithm~\cite{66} and iteratively improves the solution until reaching the maximum iteration count $\tau$. In each iteration, VNS selects a neighborhood \textit{currentN} from $N'$, generates a new solution $s$ through the shaking operation, and then performs tabu search to obtain a locally optimized solution $s^*_{TS}$. If $s^*_{TS}$ improves upon the current best solution $s^*$, it updates $s^*$ and continues searching with the same neighborhood. Otherwise, the algorithm removes \textit{currentN} from $N'$ and proceeds to the next neighborhood. When $N'$ becomes empty, it resets to the full neighborhood set $N$.

\begin{algorithm}[htbp]
\caption{Variable Neighborhood Search (VNS)}
\begin{algorithmic}[1]
\State \textit{improved} $\gets$ False
\State $N' \gets N$
\State $s^* \gets$ \textsc{CW\_Initialize}($V$)
\While{$\text{iter} \leq \tau$}
    \If{\textit{improved}}
        \State $N' \gets N'$ \Comment{Keep current neighborhood}
    \Else
        \State $N' \gets N' \setminus \{\textit{currentN}\}$ \Comment{Remove unproductive neighborhood}
        \If{$N' = \emptyset$}
            \State $N' \gets N$ \Comment{Reset neighborhoods if empty}
        \EndIf
    \EndIf
    \State $s \gets$ \textsc{Shake}(\textit{currentN}, $s^*$) \Comment{Diversification step}
    \State $(s, s_{TS}^*) \gets$ \textsc{TabuSearch}($s$, $\theta$, $\varphi$) \Comment{Intensification step}
    \If{$C(s_{TS}^*) < C(s^*)$} \Comment{Check improvement}
        \State $s^* \gets s_{TS}^*$
        \State \textit{improved} $\gets$ True
    \Else
        \State \textit{improved} $\gets$ False
    \EndIf
\EndWhile
\State \Return $s^*$
\end{algorithmic}
\end{algorithm}

\subsection{Initial Solution Construction}

Based on the C-W algorithm, this paper constructs an initial feasible solution by iteratively merging routes. The specific process is as follows:

In the first iteration, each route consists of the starting node 0, the ending node \(|N|+1\), and one customer node. When two routes are merged, one healthcare worker will visit all nodes in the new route. In each iteration, the saving value for merging routes \(r_p\) and \(r_q\) is:
\[
s_{pq} = C(r_p) + C(r_q) - C(V_r \cup V_r),
\]
where \(C(r_p)\) represents the minimum cost of route \(r_p\), including the fixed scheduling cost of a healthcare worker, travel cost, time penalty cost, and overtime cost. \(V_r \cup V_r\) denotes the union of nodes traversed by routes \(r_p\) and \(r_q\), and \(C(V_r \cup V_r)\) represents the minimum cost of the new merged route.

To simplify calculations, during the initial feasible solution merging process, only end-to-end merging of two routes is performed, without adjusting the internal node order of the merged route. For example, if \(r_p = \{0,1,2,|N|+1\}\) and \(r_q = \{0,5,4,3,|N|+1\}\), the merged new route would be \(\{0,1,2,5,4,3,|N|+1\}\) or \(\{0,5,4,3,1,2,|N|+1\}\), even if reordering the nodes could yield a lower-cost route like \(\{0,1,2,3,4,5,|N|+1\}\), it is not considered.

In each iteration, the two routes with the highest saving value are selected as candidates for merging. If the required time for the merged route does not exceed the maximum working duration \(L\) of the healthcare worker, the two routes are merged. Otherwise, the merging option with the second-highest saving value is considered, and so on. If all saving values are negative, the iteration stops, and the initialization is complete.

\subsection{Neighborhood and Perturbation}

\begin{figure}[h]
\centering
\includegraphics[width=0.7\textwidth]{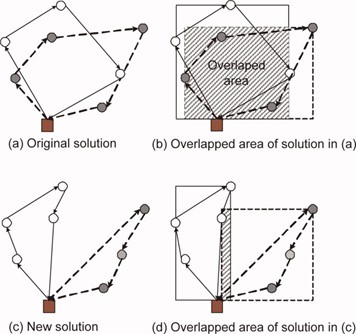}
\caption{Overlap Area}
\label{fig:overlap}
\end{figure}

During the perturbation phase, this paper constructs three neighborhoods, removing a portion of customers from the current solution and reinserting them into new routes to achieve diversified solution transitions. (i) \textit{Neighborhood N1}: Randomly select a fixed proportion (\(r\)) of customers; (ii) \textit{Neighborhood N2}: Select a fixed proportion (\(r\)) of customers with the highest relocation cost. The relocation cost of a customer refers to the reduction in total cost after removing their service point from the set~\cite{67}; (iii) \textit{Neighborhood N3}: Select customer points from routes with the largest overlapping area~\cite{68}. The region of a route is defined as the smallest rectangular area covering the starting point and all customer points on the route. As shown in Figure~\ref{fig:overlap}, the regions of routes may overlap. The overlapping area of each route is defined as the sum of its overlapping areas with all other route regions. Reducing overlapping areas may lead to solutions with shorter travel distances.

After selecting the customer points to be removed from the neighborhood, a regret-based insertion algorithm~\cite{14} is used to reinsert the selected customers into the routes. Specifically, for each customer point in the set to be reinserted, calculate its regret value, which is the difference in cost between inserting it into the best position and the second-best position. The node with the highest regret value is selected, inserted into the best position, and removed from the customer point set. This process is repeated until the customer point set is empty. The detailed implementation is described in Algorithm~\ref{algo:regret_insert}.

\begin{algorithm}[htbp]
\caption{Regret-Based Insertion}
\label{algo:regret_insert}
\begin{algorithmic}[1]
\State $V_{\text{remove}}$: Set of customers to be reinserted
\While{$N_{\text{remove}} \neq \emptyset$}
    \For{each $i \in N_{\text{remove}}$}
        \State $\text{bestP}_i \gets$ \textsc{BestInsertionCost}($i$)
        \State $\text{secondP}_i \gets$ \textsc{SecondBestInsertionCost}($i$)
    \EndFor
    \State $i^* \gets \arg \max_i (\text{secondP}_i - \text{bestP}_i)$ \Comment{Select maximum regret}
    \State \textsc{InsertCustomer}($i^*$)
    \State $N_{\text{remove}} \gets N_{\text{remove}} \setminus \{i^*\}$
\EndWhile
\end{algorithmic}
\end{algorithm}

\subsection{Local Search}
This paper employs the Tabu Search (TS) algorithm~\cite{69} to perform local search. The TS algorithm, originally proposed by Glover, can be summarized as follows: Starting from an initial solution, the algorithm iteratively moves from the current solution \( R \) to its best neighboring solution \( R' \). To prevent cycling and getting trapped in local optima, a tabu table is designed to record forbidden moves, including the tabu list and tabu tenure.

Two operators are used to generate new solution sets:
\begin{itemize}
    \item \textbf{2-Opt}: Randomly reverses a segment of a route.
    \item \textbf{Swap}: Randomly selects two routes, chooses one customer node from each, and swaps them to form two new routes.
\end{itemize}

If an edge removed in a previous step is to be reinserted, it is considered a tabu move. The tabu list \( H \) has a tabu tenure \( \theta \), which is a uniformly distributed random integer, \( \theta \sim U(\vartheta_{min}, \vartheta_{max}) \). According to the aspiration criterion, if a tabu solution yields a better objective value than the current best solution, it is allowed to override the tabu status and is used as the current solution for the next iteration, ensuring that high-quality solutions are not missed due to the tabu list. The detailed implementation is described in Algorithm~\ref{algo:local_search}.

\begin{algorithm}[htbp]
\caption{Local Search}
\begin{algorithmic}[1]
\label{algo:local_search}
\State $s'$: Initial solution
\State $s_{TS}^* \gets s$, $\text{iter} \gets 1$, $\text{stop} \gets 0$, $H \gets \emptyset$ \Comment{Initialize}
\While{$\text{iter} \leq \theta$ and $\text{stop} \leq \varphi$}
    \If{$\text{iter mod } 2 \neq 0$}
        \State $\delta(s_{TS}^*) \gets$ \textsc{2-opt}($s_{TS}^*, H$) \Comment{Alternate between 2-opt}
    \Else
        \State $\delta(s_{TS}^*) \gets$ \textsc{Swap}($s_{TS}^*, H$) \Comment{and swap moves}
    \EndIf
    \State $s' \gets$ \textsc{SelectBestSolution}($\delta(s_{TS}^*)$)
    \State \textsc{Update}($H$) \Comment{Update tabu list}
    \If{$C(s') < C(s_{TS}^*)$}
        \State $s_{TS}^* \gets s'$
        \State $\text{stop} \gets 0$ \Comment{Reset non-improvement counter}
    \Else
        \State $\text{stop} \gets \text{stop} + 1$
    \EndIf
    \State $\text{iter} \gets \text{iter} + 1$
\EndWhile
\State \Return $s, s_{TS}^*$
\end{algorithmic}
\end{algorithm}

\section{Numerical Experiments}
\label{sec:exps}

\subsection{Experiment Setup}
To the best of our knowledge, there are no benchmark studies or algorithms specifically designed for the VRASP problem in the existing literature. Therefore, this study employs both CPLEX 12.10.0 to obtain exact solutions and the VNS algorithm to derive approximate solutions, comparing the results of the two approaches. Additionally, both methods are applied to solve the deterministic VRASP (\(P_0\)) and the stochastic VRASP (\(P_1\)) under SAA modeling, with the results of the two models analyzed and compared. All experiments were conducted on a personal computer equipped with a 2.1 GHz AMD processor and 16 GB of RAM.

For the numerical experiments, test instances were constructed based on parameter settings from recent related literature. The starting node is located at (0, 0), and customer points are uniformly distributed within the \([0, 50]^2\) area. According to NAHC3 statistics, the average daily number of patients served by a healthcare worker in home healthcare is approximately 6. Thus, the number of healthcare workers is set to \(|K| = \lceil |N| / 6 \rceil\). The maximum working duration for each healthcare worker is \(L = 480\) minutes, with a fixed scheduling cost \(c' = 100\), a unit overtime cost \(c'' = 1\), and a unit tardiness cost \(c''' = 2\). The travel cost \(t_{ij}\) between customers \(i\) and \(j\) is proportional to the Euclidean distance, with a scaling factor \(l = 0.5\). The travel time \(t_{ij}\) is also distance-dependent, with the scaling factor following a uniform distribution \(\gamma \sim U(0.5, 1.5)\). The service time \(ts_i\) at customer \(i\) follows a truncated log-normal distribution over the interval \([30, 90]\), with \(\sigma = 0.5\mu\).

The experiments include five problem scales with \(|N| = 5, 10, 20, 30, 40\). For each scale, 10 test instances are generated. For each instance, random travel and service times are sampled \(|M| = 30, 50, 80, 100\) times to solve the SAA model. For the deterministic VRASP, the sample means of the random travel times \(t_{ij}\) and service times \(ts_i\) are substituted into \(P_0\) for solving.

The VNS algorithm is configured with a maximum iteration count \(t = 1000\). During the shaking phase, the proportion of customers selected from neighborhoods \(N1\) and \(N2\) is \(r = 0.2\). The local search process has two termination criteria: (1) a maximum iteration limit \(\theta = 1000\), and (2) a maximum iteration limit \(\varphi = 5 \times |N|\) if no improvement in the optimal solution is observed.

\subsection{Running Time Comparison}\label{ssec:run_time}

This study employs both CPLEX 12.10.0 and the VNS algorithm to solve the deterministic VRASP (\(P_0\)) and the stochastic VRASP (\(P_1\)). For \(P_1\), the sample size \(|M|\) in the SAA method significantly impacts the computational time. Thus, four sample sizes \(|M| = 30, 50, 80, 100\) are tested. In \(P_0\), the random travel times \(t_{ij}\) and service times \(ts_i\) are replaced by their sample means. Since the computation time for the mean is not significantly affected by the sample size, the largest sample size \(|M| = 100\) is selected to ensure the mean values closely approximate the true values under the law of large numbers. Five problem scales \(|N| = 5, 10, 20, 30, 40\) are tested, with 10 instances generated for each scale. The average computational time for each instance is recorded, as shown in Table~\ref{tab:runtime}.

Comparing the computational times of \(P_0\) and \(P_1\), it is evident that \(P_0\) is solved significantly faster than \(P_1\) by both CPLEX and VNS. This is because the SAA method requires evaluating constraints and objectives for each sample, leading to longer computation times as the sample size increases. In contrast, the deterministic model \(P_0\) considers only one sample (the mean of all samples), drastically reducing the computational effort.

When comparing CPLEX and VNS, CPLEX outperforms VNS for small-scale instances (e.g., \(|N|=5, |K|=1\)), highlighting its efficiency for such problems. However, the heuristic algorithm in this study sets high iteration limits for both the outer VNS loop (\(\tau = 1000\)) and the inner local search (\(\theta = 1000\)), which inherently limits further reduction in runtime. CPLEX's computation time increases superlinearly with problem size. For \(|N| = 10\), CPLEX-P0 remains faster than VNS-P0 (69.75 vs. 2060.58 seconds), but for \(|N| = 20\), it becomes significantly slower (36561.85 vs. 5223.19 seconds). For \(|N| = 30\), CPLEX fails to produce results within a reasonable time. The sample size also greatly affects CPLEX's performance. For \(|N| = 10\) and \(|M| = 30\), CPLEX-P1 is already slower than VNS-P1 (81258.94 vs. 2060.58 seconds), and for \(|M| = 50\), CPLEX cannot complete the computation. These results demonstrate the clear advantage of the proposed VNS algorithm in solving medium- to large-scale problems efficiently.

\begin{table}[h]
\renewcommand{\arraystretch}{1.2} 
\setlength{\tabcolsep}{8pt} 
\centering
\caption{Comparison of computational time (seconds). ``-" indicates no solution obtained within 30 hours (108,000 seconds)}
\label{tab:runtime}
\begin{tabular}{lcccccc}
\hline
Model & Sample Size & \multicolumn{5}{c}{Problem Size (Number of Customers)} \\ 
 & & 5 & 10 & 20 & 30 & 40 \\ \hline
CPLEX-P0 & - & 0.06 & 69.75 & 36561.85 & - & - \\ 
VNS-P0 & - & 429.44 & 527.97 & 638.08 & 864.76 & 1100.14 \\ \hline
CPLEX-P1 & 30 & 9.16 & 81258.94 & - & - & - \\ 
 & 50 & 16.09 & - & - & - & - \\ 
 & 80 & 17.68 & - & - & - & - \\ 
 & 100 & 21.55 & - & - & - & - \\ \hline
VNS-P1 & 30 & 927.91 & 2060.58 & 5223.19 & 10490.09 & 23518.10 \\ 
 & 50 & 1550.42 & 3864.23 & 7457.93 & 16965.71 & 36742.74 \\ 
 & 80 & 2671.89 & 6652.06 & 12596.42 & 27551.95 & 52544.61 \\ 
 & 100 & 4158.57 & 11533.80 & 20775.14 & 38538.23 & 71578.35 \\ \hline
\end{tabular}
\end{table}

\subsection{Deterministic vs. Stochastic Models}

Based on the results from Section~\ref{ssec:run_time}, the computational time of CPLEX increases dramatically with both problem size and sample size. Therefore, this study selects \(|N| = 10\) and \(|M| = 30\) as the standard test case. The cost is calculated as the average total cost across all samples for given service routes \(x_{ijk}\) and service start times \(s_i\). The gap between the two models is computed as:
\[
\text{gap} = \frac{P_0 - P_1}{P_0}
\]

Tables~\ref{tab:1} and~\ref{tab:2} compare the results of solving \(P_0\) and \(P_1\) using CPLEX and VNS, respectively. In both methods, the total cost obtained from \(P_1\) is consistently lower than that from \(P_0\). One possible explanation is that the SAA method accounts for sample diversity by minimizing the sum of total costs across all samples, while the deterministic model simplifies the problem by using sample averages for decision-making. Although averaging simplifies the problem, the optimal solution based on averaged values may not minimize the cost for individual samples, leading to higher overall costs.

\begin{table}[h]
\renewcommand{\arraystretch}{1.2} 
\setlength{\tabcolsep}{12pt} 
\centering
\caption{CPLEX-P0 vs. CPLEX-P1 Performance Comparison}
\label{tab:1}
\begin{tabular}{cccccc}
\hline
Instance & \multicolumn{2}{c}{CPLEX-P0} & \multicolumn{3}{c}{CPLEX-P1} \\
\cline{2-6}
 & Cost & Time (s) & Cost & Gap & Time (s) \\
\hline
1 & 322.41 & 466.26 & 309.43 & 4.02\% & 2401.01 \\
2 & 290.14 & 477.71 & 283.00 & 2.46\% & 2111.83 \\
3 & 331.68 & 474.09 & 328.67 & 0.91\% & 1977.22 \\
4 & 326.66 & 442.45 & 315.23 & 3.50\% & 1967.78 \\
5 & 298.47 & 440.91 & 287.30 & 3.74\% & 1988.16 \\
6 & 306.02 & 405.68 & 301.53 & 1.47\% & 2011.05 \\
7 & 314.35 & 442.84 & 313.57 & 0.25\% & 1972.20 \\
8 & 324.00 & 443.25 & 311.90 & 3.73\% & 1964.53 \\
9 & 319.77 & 399.86 & 316.77 & 0.94\% & 1831.03 \\
10 & 320.74 & 436.65 & 319.37 & 0.43\% & 2381.01 \\
\hline
Avg. & 315.42 & 442.97 & 309.68 & 1.82\% & 2060.58 \\
\hline
\end{tabular}
\end{table}

\begin{table}[h]
\renewcommand{\arraystretch}{1.2} 
\setlength{\tabcolsep}{12pt} 
\centering
\caption{VNS-P0 vs. VNS-P1 Performance Comparison}
\label{tab:2}
\begin{tabular}{cccccc}
\hline
Instance & \multicolumn{2}{c}{VNS-P0} & \multicolumn{3}{c}{VNS-P1} \\
\cline{2-6}
 & Cost & Time (s) & Cost & Gap & Time (s) \\
\hline
1 & 321.00 & 132.73 & 308.00 & 4.05\% & 94378.95 \\
2 & 290.14 & 69.20 & 280.14 & 3.45\% & 73957.44 \\
3 & 330.00 & 35.03 & 326.00 & 1.21\% & 80921.37 \\
4 & 325.66 & 18.20 & 313.66 & 3.68\% & 88489.21 \\
5 & 298.31 & 47.38 & 286.31 & 4.02\% & 69580.46 \\
6 & 304.64 & 136.85 & 300.64 & 1.31\% & 72581.36 \\
7 & 312.00 & 48.90 & 312.00 & 0.00\% & 90002.71 \\
8 & 321.96 & 10.59 & 310.96 & 3.42\% & 85439.67 \\
9 & 318.00 & 330.37 & 315.00 & 0.94\% & 77584.31 \\
10 & 320.73 & 55.65 & 316.73 & 1.25\% & 79653.92 \\
\hline
Avg. & 314.24 & 88.49 & 306.94 & 2.32\% & 81258.94 \\
\hline
\end{tabular}
\end{table}

\subsection{Exact Methods vs. Heuristic Algorithms}

Tables~\ref{tab:3} and~\ref{tab:4} compare the results obtained from CPLEX and VNS. The gap between the two methods is calculated as:
\[
\text{gap} = \frac{\text{VNS} - \text{CPLEX}}{\text{VNS}}
\]

The results show that while CPLEX slightly outperforms VNS for both \(P_0\) and \(P_1\), the difference is minimal. Across 10 test instances, the average gaps are only 0.37\% and 0.56\%, respectively, demonstrating the accuracy of the VNS algorithm. Since VRASP is an NP-hard problem, exact methods like CPLEX are only feasible for small-scale instances. In practice, it is more reasonable to balance solution quality and computational efficiency by using heuristic algorithms, which provide near-optimal solutions with acceptable accuracy.

\begin{table}[h]
\renewcommand{\arraystretch}{1.2} 
\setlength{\tabcolsep}{12pt} 
\centering
\caption{CPLEX-P0 vs. CPLEX-P1 Performance Comparison}
\label{tab:3}
\begin{tabular}{cccccc}
\hline
Instance & \multicolumn{2}{c}{CPLEX-P0} & \multicolumn{3}{c}{CPLEX-P1} \\
\cline{2-6}
 & Cost & Time (s) & Cost & Gap & Time (s) \\
\hline
1 & 321.00 & 132.73 & 322.41 & 0.44\% & 466.26 \\
2 & 290.14 & 69.20 & 290.14 & 0.00\% & 477.71 \\
3 & 330.00 & 35.03 & 331.68 & 0.51\% & 474.09 \\
4 & 325.66 & 18.20 & 326.66 & 0.31\% & 442.45 \\
5 & 298.31 & 47.38 & 298.47 & 0.05\% & 440.91 \\
6 & 304.64 & 136.85 & 306.02 & 0.45\% & 405.68 \\
7 & 312.00 & 48.90 & 314.35 & 0.75\% & 442.84 \\
8 & 321.96 & 10.59 & 324.00 & 0.63\% & 443.25 \\
9 & 318.00 & 330.37 & 319.77 & 0.55\% & 399.86 \\
10 & 320.73 & 55.65 & 320.74 & 0.00\% & 436.65 \\
\hline
Avg. & 314.24 & 88.49 & 315.42 & 0.37\% & 442.97 \\
\hline
\end{tabular}
\end{table}

\begin{table}[h]
\renewcommand{\arraystretch}{1.2} 
\setlength{\tabcolsep}{12pt} 
\centering
\caption{CPLEX-P0 vs. CPLEX-P1 Performance Comparison}
\label{tab:4}
\begin{tabular}{cccccc}
\hline
Instance & \multicolumn{2}{c}{CPLEX-P0} & \multicolumn{3}{c}{CPLEX-P1} \\
\cline{2-6}
 & Cost & Time (s) & Cost & Gap & Time (s) \\
\hline
1 & 308.00 & 94378.95 & 309.43 & 0.46\% & 2401.01 \\
2 & 280.14 & 73957.44 & 283.00 & 1.01\% & 2111.83 \\
3 & 326.00 & 80921.37 & 328.67 & 0.81\% & 1977.22 \\
4 & 313.66 & 88489.21 & 315.23 & 0.50\% & 1967.78 \\
5 & 286.31 & 69580.46 & 287.30 & 0.34\% & 1988.16 \\
6 & 300.64 & 72581.36 & 301.53 & 0.30\% & 2011.05 \\
7 & 312.00 & 90002.71 & 313.57 & 0.50\% & 1972.20 \\
8 & 310.96 & 85439.67 & 311.90 & 0.30\% & 1964.53 \\
9 & 315.00 & 77584.31 & 316.77 & 0.56\% & 1831.03 \\
10 & 316.73 & 79653.92 & 319.37 & 0.83\% & 2381.01 \\
\hline
Avg. & 306.94 & 81258.94 & 308.68 & 0.56\% & 2060.58 \\
\hline
\end{tabular}
\end{table}

\section{Conclusion}\label{sec:conclusion}
This paper introduces a novel \textit{Vehicle Routing and Appointment Scheduling Problem} (VRASP) that simultaneously optimizes service routes and scheduling in home healthcare systems. We first develop a deterministic VRASP model minimizing total costs including scheduling, travel, lateness penalties, and overtime work. The model is then extended to account for travel and service time uncertainties using the Sample Average Approximation (SAA) method. To solve these models, we design a Variable Neighborhood Search (VNS) heuristic and validate its accuracy and efficiency through extensive computational experiments comparing against CPLEX exact solutions.

The theoretical contributions are threefold: (i) VRASP represents a new joint optimization variant that enriches vehicle routing problem research; (ii) The proposed VNS-based heuristic demonstrates both solution accuracy for small instances and computational efficiency for medium-to-large instances; (iii) The methodology establishes a foundation for solving related optimization problems. Practically, our approach addresses critical challenges in home healthcare by accounting for the distinct needs of service providers and recipients, offering a flexible solution that reduces unnecessary waiting through joint scheduling and routing optimization under uncertainty.

While this study advances VRP research, several limitations suggest promising future directions. First, rather than fixed appointment times, future work could consider time windows for service initiation. Second, alternative methods like stochastic dynamic programming or chance-constrained programming could complement our SAA approach for handling uncertainties. Third, though VNS shows strong agreement with CPLEX for small instances, its solution quality for large-scale problems requires further validation due to the lack of comparable benchmarks. These extensions would further enhance the model's practical applicability in home healthcare logistics.

\newpage
\bibliographystyle{plain}
\bibliography{main}

\end{document}